\documentclass[a4,reqno,11pt]{amsart}

\usepackage{amsmath, amsthm, amssymb, amscd, mathabx}
\usepackage{mathdots, mathrsfs}

\usepackage{graphicx, verbatim, url, enumerate, relsize}
\usepackage[dvipsnames]{xcolor}
\usepackage[normalem]{ulem} 

\usepackage{tikz}
\usepackage{tikz-cd}
\usetikzlibrary{matrix}

\usepackage[pdftex]{hyperref}
\usepackage{cleveref}
\renewcommand{\setminus}{\smallsetminus}
\newcommand{\ovl}{\overline}

\newcommand{\wh}[1]{\widehat{#1}}

\newcommand{\la}{\langle}
\newcommand{\ra}{\rangle}

\newcommand{\qq}{\mathbb Q}

\newcommand{\cc}{\mathbb C}

\newcommand{\nat}{\mathbb N}
\newcommand{\zz}{\mathbb Z}
\newcommand{\rr}{\mathbb R}

\newcommand{\mc}[1]{\mathcal{#1}}
\newcommand{\mf}[1]{\mathfrak{#1}}
\newcommand{\mfp}{\mathfrak{p}}
\newcommand{\mfm}{\mathfrak{m}}

\newcommand{\mg}[1]{{#1}^{\times}}
\newcommand{\sq}[1]{{#1}^{\times 2}}

\newcommand{\vf}{\varphi}

\renewcommand{\geq}{\geqslant}

\DeclareMathOperator{\supp}{\mathsf{Supp}}

\DeclareMathOperator{\Spec}{\mathsf{Spec}}

\renewcommand{\dim}{\operatorname{\mathsf{dim}}}

\numberwithin{equation}{section}
\newtheorem{thm}[equation]{Theorem}
\newtheorem{prop}[equation]{Proposition}
\newtheorem{cor}[equation]{Corollary}

\newtheorem{lem}[equation]{Lemma}

\newtheorem*{theorem*}{Theorem}
\newtheorem*{qu*}{Question}

\theoremstyle{definition}

\newtheorem{exs}[equation]{Examples}

\swapnumbers

\title[Quadratic forms over semiglobal fields]{A local-local principle for quadratic forms over semiglobal fields}

\author{Karim Johannes Becher}
\author{David Grimm}
\author{Parul Gupta}

\email{karimjohannes.becher@uantwerpen.be}
\email{david.grimm@usach.cl}
\email{parul.gupta1@snu.edu.in}

\address{University of Antwerp, Department of Mathematics, Middelheim\-laan~1, 2020 Antwerpen, Belgium.}

\address{Universidad de Santiago de Chile, Facultad de Ciencias, Avenida Libertador Bernardo O'Higgins nº 3363, Estaci\'on Central, Santiago, Chile}

\address{Department of Mathematics, Indian Institute of Technology Indore, Simrol, Khandwa Road, Indore 453 552 India}

\date{16.11.2026}


\begin{document}

\begin{abstract}\noindent
Given a quadratic form over a semiglobal field, the set of discrete valuations on the field for which the form is anisotropic over the completion is either empty, or it contains a valuation extending that of the base field with a transcendental residue field extension.
In dimension at least $3$, this gives a refinement of the  local-global principle for isotropy.
In dimension $2$, the result yields a description of the set of elements that are locally squares everywhere.
\end{abstract}

\maketitle

\section{Introduction}

By a \textit{local-global principle} we mean a result that deduces a property of some algebraic structure $\Phi$ (e.g.~a quadratic form) defined over a field $F$ from the hypothesis that this property holds for $\Phi$ over every member of some distinguished system of field extensions of $F$.
The prototype of a local-global principle concerning quadratic forms is the Hasse--Minkowski Theorem.
Before recalling it, we fix some terminology.

For a field $F$, we denote by $\Omega(F)$ the 
set of valuations on $F$ with value group $\zz$, which we call \textit{$\zz$-valuations}.
For $v\in\Omega(F)$, we denote by $F^v$ the completion of $F$ with respect to $v$.

By a \textit{function field in one variable} we mean a finitely generated field extension of transcendence degree $1$.
A \textit{global field} is a field which is a finite extension of $\qq$ or a function field in one variable over a finite field. 

A quadratic form over a field is called \textit{isotropic} if it represents zero nontrivially, and   \textit{anisotropic} otherwise.
A quadratic form $\vf$ over the field $F$ is \textit{totally indefinite} if it represents positive and negative elements with respect to any field ordering of $F$; this condition is trivial if the field $F$ has no orderings.

\begin{thm}[Hasse--Minkowski]\label{HM}
    Let $\vf$ be a quadratic form over a global  field $F$.
    Then $\vf$ is isotropic over $F$ if and only if it is totally indefinite and $\vf$ is isotropic over $F^v$ for every $v\in\Omega(F)$.
\end{thm} 
Here, the set $\Omega(F)$ cannot be replaced by any proper subset:
For any $v_0\in\Omega(F)$, one can construct a $4$-dimensional quadratic form over $F$ which is anisotropic over $F^{v_0}$ but isotropic over $F^v$ for all $v\in\Omega(F)\setminus\{v_0\}$.
For example,  the indefinite quadratic form $X_1^2+X_2^2+X_3^2-7X_4^2$ over $\qq$ is anisotropic over $\qq_2$ but isotropic over $\qq_p$ for every odd prime number $p$.

There is an analogue to the Hasse-Minkowski Theorem for function fields in one variable over $\rr$, due to Ernst Witt \cite[Satz 22]{Wit37}. 
The following version is covered by \cite[Theorem I]{ELP73}
together with \cite[Theorem 9.4]{Kne76}.
(It holds trivially if the field $F$ has no  orderings.)

\begin{thm}[Witt]\label{Witt}
Let $F/K$ be a function field in one variable where $K$ is real closed. 
A quadratic form in at least three variables over $F$ is isotropic if and only if it is totally indefinite, if and only if it is isotropic over $F^v$ for all but finitely many $v\in\Omega(F)$.
\end{thm}
Trivially, this local-global principle remains valid when $\Omega(F)$ is replaced by any cofinite subset.
Another difference to \Cref{HM} is that binary quadratic forms are not covered. 
In fact, for $t\in F$ transcendental over $K$, the totally indefinite quadratic form $X_1^2-(n+t^2)X_2^2$ is 
anisotropic over such a field $F$ for almost all  $n\in\nat\setminus\{0\}$.

We now discuss two types of fields where local-global principles for quadratic forms (and other principal homogeneous spaces)  were established recently, in \cite[Theorem 3.1]{CTPS12} and \cite[Corollary 4.7]{HHK15}, respectively.

\begin{thm}\label{CTPS-HHK}
In either of the following situations, a quadratic form $\vf$ over $F$ in at least three variables is isotropic if and only if it is isotropic over $F^v$ for every $v\in\Omega(F)$.
\begin{enumerate}[$(i)$]
    \item $F/K$ is a function field in one variable where $K$ is the fraction field of a complete discrete valuation ring $T$ with $2\in\mg{T}$. \\ \phantom{x}\hfill {\rm (Colliot-Thélène--Parimala--Suresh)}
    \item $F$ is the fraction field of a regular or equicharacteristic $2$-dimensional complete local ring $T$ with $2\in\mg{T}$. \hfill{\rm (Harbater--Hartmann--Krashen)}
\end{enumerate}
\end{thm}
The proofs rely on field patching techniques, introduced in \cite{HH10} and \cite{HHK09}.
A key point in the approach is a geometric interpretation of $F$ as the function field of an arithmetic surface over $\Spec(T)$.
Furthermore, the flexibility in the choice of this surface is crucial.
This relies on the embedded resolution of singularities for surfaces due to  Lipman.

In contrast to the situations of \Cref{HM} and \Cref{Witt}, the set $\Omega(F)$ in the two scenarios of \Cref{CTPS-HHK} is vast and lacks a simple geometric parametrization.
For example, for any $a\in\mg{F}$, the set $\{v\in\Omega(F)\mid v(a)\neq 0\}$ is finite in \Cref{HM} and \Cref{Witt}, whereas it is infinite in \Cref{CTPS-HHK}.
It is therefore natural to ask for a good choice of a subset of $\Omega(F)$ which is sufficient for the local-global principle.
Our aim is to provide such a refinement. 
The first possibility coming to mind in the context of situation $(i)$ may be to restrict to valuations $v\in\Omega(F)$ which are trivial on $K$.
This, however, is generally not possible, even not for $K=\cc(\!(t)\!)$ and $F=K(X)$; see ~\cite{Gup21}.

The study of the two scenarios in \Cref{CTPS-HHK} has inspired a new term for the corresponding fields $F$.
Let $T$ be an integrally closed, noetherian, henselian local subring of $F$.
Assume that one of the following holds:
\begin{enumerate}[$(i)$] 
\item $T$ is $1$-dimensional, hence a discrete valuation ring, and $F$ is a function field in one variable over the fraction field of $T$.
\item $T$ is $2$-dimensional and $F$ is the fraction field of $T$.
\end{enumerate}
In this situation we call $F/T$ a \textit{semiglobal field of type $(i)$}, respectively \textit{of type $(ii)$}.

For an arbitrary valuation $v$ on $F$, we denote by $\mc{O}_v$ the corresponding valuation ring, by $\mfm_v$ its maximal ideal and by $Fv$ its residue field $\mc{O}_v/\mfm_v$.
We refer to \cite[Sect.~2.4]{EP05} for the notion of completion of a valued field. 
For a valuation $v$ on $F$, we denote by $F^v$ the corresponding completion of $F$ with respect to $v$ and by $\hat{v}$ the natural extension of $v$ to $F^v$.
If $v$ is a discrete valuation on $F$, then the valuation ring $\mc{O}_{\hat{v}}$ is the $\mfm_v$-adic completion of $\mc{O}_v$.

Recall that $T$ is a local subring of $F$.
Let $\mfm$ denote its maximal ideal and $k=T/\mfm$.
We denote by $\Omega_T(F)$ the subset of $\Omega(F)$ given by $\zz$-valuations $v$ for which $T\subseteq \mc{O}_v$ and $\mfm_v\cap T=\mfm$. 
Note that, for $v\in\Omega_T(F)$, the residue field $Fv$ is a an extension of $k$.
Let 
\[\Omega_T^\circ(F)=\{v\in\Omega_T(F)\mid Fv/k\mbox{ is a function field in one variable}\}\,.\]
We aim to replace $\Omega(F)$ by  $\Omega_T^\circ(F)$ in \Cref{CTPS-HHK}.
This will be achieved in \Cref{2-foldlgpi} by combining \Cref{CTPS-HHK} with the following local-local principle, which is at the center of this article:

\begin{thm}[{\Cref{1-2-dim}}]\label{main}
Let $F/T$ be a semiglobal field where $T$ is excellent and with $2\in\mg{T}$.
Let $\vf$ be a quadratic form over $F$.
Assume that $\vf$ is isotropic over $F^v$ for every $v\in\Omega_T^\circ(F)$.
Then $\vf$ is isotropic over $F^v$ for every rank-$1$ valuation $v$ on $F$ with $T\subseteq \mc{O}_v$.
\end{thm} 

We refer to \Cref{main} as a \textit{local-local principle} because it deduces a local property (isotropy over all completions) from a heavily restricted subset of local conditions.
By a \emph{rank-$1$ valuation} we mean a nontrivial valuation taking values in $\rr$.

While the hypothesis of \Cref{main} only requires testing isotropy over completions with respect to specific discrete valuations (those in $\Omega_T^\circ(F)$), the conclusion holds for all rank-$1$ valuations that are non-negative on $T$. Although this inclusion of non-discrete valuations is primarily a byproduct of our methodology, it yields a strictly stronger statement.
For the refinement in \Cref{CTPS-HHK}, we only need the conclusion for the discrete valuations in $\Omega_T(F)$.

The hypothesis in \Cref{main} that $T$ is excellent is satisfied if, for example, $T$ is complete or of characteristic zero; see~\Cref{EX:excellent}.
For the refinement in \Cref{CTPS-HHK}, it would be sufficient to state the conclusion in \Cref{main} for all discrete valuations $v\in\Omega(F)$.
The version stated has an additional benefit in the case $(i)$: it shows that the local-global principle from \cite[Theorem 3.1]{CTPS12} can be derived from (and hence seen as a special case of) Mehmeti's local-global principle  \cite[Cor. 3.19 \& Rem. 3.20]{Meh19}.

Our arguments leading to the proof of \Cref{main} are based on classical commutative algebra and algebraic geometry, including Lipman's  resolution of singularities for normal surfaces as the only deep ingredient.

\section{Forms over two-dimensional regular local rings}\label{S:prelim} 

For quadratic forms defined over a two-dimensional regular local ring, we study their isotropy behavior over the fraction fields of the completions with respect to the height-one prime ideals.

Every regular local ring is  a domain; see \cite[Cor.~10.14]{Eisenbud}.
Let $A$ be a local ring and let $\mfm$ be its maximal ideal.
We call $A/\mfm$ the \textit{residue field of $A$}.
We denote by $\wh{A}$ its completion with respect to $\mfm$, which is again a local ring with the same residue field.
If $A$ is regular, then so is $\wh{A}$; see \cite[Cor.~10.12]{Eisenbud}.

Let $R$ be a domain and let $F$ be its fraction field. Let $\mf p$ be a prime ideal of $R$. 
We denote by $R_{\mf p}$ the localization of $R$ with respect to $\mf p$ and by $\kappa_\mf p$ the residue field of $R_{\mf p}$, which is given by the fraction field of $R/\mfp$. The residue field of $\wh{R}_\mfp$ is again $\kappa_\mfp$.
When $\wh{R}_\mathfrak{p}$ is a domain (e.g.~when~$R$ is regular), we denote by $F^\mf p$  the fraction field of $\wh{R}_{\mf p}$. 

\begin{prop}[Krull-Akizuki]\label{KA}
Let $A$ be a noetherian local domain of Krull dimension $1$. 
Let $K$ be the fraction field of $A$, $L/K$ be a finite field extension and $B$ the integral closure of $A$ in $L$.
Then $B$ is a Dedekind domain that is not a field. Moreover, if $A$ is henselian, then $B$ is a henselian discrete valuation ring.
\end{prop}
\begin{proof}
Since $A$ is noetherian of Krull dimension $1$, so is $B$, by \cite[Chap.~VII, \S2, n$^\circ$5, Prop.~5]{BOU-AC5to7}.
Furthermore, since $B$ is integrally closed, it follows that $B$ is a Dedekind domain which is not a field.

Suppose now that $A$ is henselian.
Then \cite[Cor.~30.5]{Nag75} yields that $B$ is local, and hence a discrete valuation ring.
Furthermore, it follows by \cite[Cor.~43.13]{Nag75} that $B$ is henselian.
\end{proof}

Given $n\in\nat$ and $a_1,\dots,a_n\in\mg{F}$, we denote by $\la a_1,\dots,a_n\ra_F$ (or simply $\la a_1,\dots,a_n\ra$) the $n$-ary quadratic form $a_1X_1^2+\ldots+a_nX_n^2$ over $F$.
We say that $\la a_1,\dots,a_n\ra_F$ is \textit{isotropic} if this form represents zero non-trivially over $F$, and \textit{anisotropic} otherwise.

\begin{lem}\label{unimodulargeneral}
Let $R$ be a $2$-dimensional regular local ring with maximal ideal $\mf m$ and $2\in\mg{R}$. Let $F$ denote the fraction field and $\kappa$ the residue field of $R$.
 Let $a_1, \ldots, a_r\in\mg R$ be such that $\la \ovl{a}_1, \ldots, \ovl{a}_r\ra_{\kappa}$ is isotropic.
Let $v$ be a rank-$1$ valuation on $F$ with $R\subseteq \mc{O}_v$. 
If $\mfm_v \cap R\neq \mfm$, then assume that the residue field $Fv$ carries a henselian $\zz$-valuation.
Then $\la a_1, \ldots, a_r\ra_{F^v}$ is isotropic.
\end{lem}

\begin{proof}
Since $v$ is a rank-$1$ valuation, the completion $F^v$ is henselian with respect to the natural extension of $v$, by \cite[Theorem 1.3.1]{EP05}.

Set $\mfp=\mfm_v\cap R$.
If $\mfp=\mfm$, then $\kappa$ embeds into $Fv$, whereby the hypothesis yields that $\la \ovl{a}_1, \ldots, \ovl{a}_r\ra_{Fv}$ is isotropic, and since $F^v$ is henselian with respect to the natural extension of $v$, we conclude that $\la a_1,\dots,a_r\ra_{F^v}$ is isotropic.

Assume now that $\mfp\neq \mfm$.
Since $v$ is a nontrivial valuation, we have that $\mfp\neq \{0\}$.
Hence $\mfp$ is a height-one prime ideal of $R$.
It follows that $R_{\mfp}$ is a discrete valuation ring of $F$ contained in $\mc{O}_v$.
In particular, $R_\mfp$ is a maximal proper subring of $F$.
Since $R_\mfp\subseteq\mc{O}_v\subsetneq F$, we conclude that $\mc{O}_v=R_\mfp$.
In particular, $Fv=\kappa_{\mfp}$. 
Let $A$ denote the normalization of $R/\mf p$.
Then, by \Cref{KA}, $A$ is a Dedekind domain that is not a field.
By the hypothesis, $Fv$ carries a henselian $\zz$-valuation $w$.
By \cite[Prop.~2.2]{BGVG14}, then $w$ is the unique $\zz$-valuation on $Fv$.
Since the localization of $A$ at any maximal ideal is a discrete valuation ring of $Fv$, it follows that $A$ is  local and thus a discrete valuation ring of $Fv$.
    Hence $A=\mc{O}_w$, so $A$ is henselian.

The residue field of $A$ contains $\kappa$, the residue field of $R$. Since $\la \ovl{a}_1, \ldots, \ovl{a}_r\ra_{\kappa}$ is isotropic, $A$ is henselian and $2\in\mg{A}$, we get that $\la \ovl{a}_1, \ldots, \ovl{a}_r\ra_{\kappa_\mfp}$ is isotropic. 
Since $\kappa_\mfp=Fv$ and since $F^v$ is henselian with respect to the natural extension of $v$, we conclude
that $\la a_1,\dots,a_r\ra_{F^v}$ is isotropic.
\end{proof}

Assume that $R$ is a noetherian unique factorization domain.
Let $\mc P_R$ be the set of height-one prime ideals of $R$. 
For every $\mf p \in \mc P_R$,  the ring $R_\mf p$ is a discrete valuation ring, and we then denote by $v_{\mfp}$ the corresponding $\zz$-valuation on $F$.
To any $a\in \mg {F}$, we associate its \textit{support} 
$$\supp_R(a) = \{ \mf{ p} \in \mc P_R \mid  v_{\mf p}(a) \neq 0\}\,. $$
Any regular local ring is a unique factorization domain, by \cite[Theorem 19.19]{Eisenbud}.

Let $\varphi$ be a regular quadratic form over $F$. We say that \textit{$\varphi$ is unimodular over $R$} if there exists $r\in\nat$ and $a_1, \ldots, a_r\in \mg R$ such that $ \varphi \simeq \la a_1, \ldots, a_r\ra$. Moreover, if $R$ is local with maximal ideal $\mf m$ and  residue field $\kappa$ and if $\varphi = \la a_1, \ldots, a_r\ra$ with $a_1, \ldots, a_r \in \mg R$, then we call the diagonal form $\la \ovl{a}_1, \ldots, \ovl{a}_r\ra$  over $\kappa$ the \textit{residue form of $\varphi$ with respect to $\mf m$} and we denote it by $\ovl{\varphi}$. 
Note that we define residue forms only for diagonal forms, and a priori not up to isometry.

\begin{lem}\label{residueisotropic}
Let $R$ be a $2$-dimensional regular local ring with maximal ideal $\mf m$, fraction field $F$, residue field $\kappa$ and $2\in\mg{R}$. Let $(\pi, \delta)$ be a parameter system of $R$. Let $\varphi_0,\varphi_1,\varphi_2,\varphi_3 $ be unimodular diagonal quadratic forms over $R$ and let 
\[\varphi= \varphi_0 \perp \pi \varphi_1 \perp \delta \varphi_2 \perp \pi\delta \varphi_3.\]
Let $F^{(\delta)}$ denote the fraction field of the completion of $R$ with respect to the height-$1$ prime ideal $(\delta)$.
Assume that $\varphi$ is isotropic over $F^{(\delta)}$. 
Then the following hold:
\begin{enumerate}[(a)]
    \item Over $\kappa$, at least one of the residue forms modulo $\mathfrak{m}$ of $\varphi_0,{\varphi}_1,{\varphi}_2,{\varphi}_3$ is isotropic.
    
   \item 
   For any rank-$1$ valuation $v$ on $F$  such that $R\subseteq \mc{O}_v$ and, either $\mf m_v\cap R=\mf m$, or the residue field $Fv$ carries a henselian  discrete valuation, $\varphi$ is isotropic over $F^v$.
    \end{enumerate}
\end{lem}

\begin{proof} 
 Since $\varphi$ is isotropic over $F^{(\delta)}$, it follows by \cite[Lemma 19.5]{EKM08} that one of the residue forms $\ovl{\varphi_0} \perp \ovl{\pi}\cdot \ovl{\varphi_1}$ and $\ovl{\varphi_2}\perp \ovl{\pi} \cdot \ovl{\varphi_3}$ modulo $(\delta)$ is isotropic over $\kappa_{(\delta)}$. 
 Since $\delta$ is a parameter of $R$,  $R/(\delta)$ is a discrete valuation ring with fraction field $\kappa_{(\delta)}$,  maximal ideal $(\ovl{\pi}) = (\pi+R\delta)$ and residue field $\kappa=R/(\pi,\delta)$. 
Over
$R/(\delta)$, we obtain
 diagonal unimodular forms $\ovl{\varphi_0},\ldots, \ovl{\varphi_3}$.
Since one of the forms $\ovl{\varphi_0} \perp \ovl{\pi}\cdot \ovl{\varphi_1}$ and $\ovl{\varphi_2}\perp \ovl{\pi} \cdot \ovl{\varphi_3}$ is isotropic over $\kappa_{(\delta)}$ and since $\ovl{\pi}$ is a uniformizer of $R/(\delta)$, by \cite[Lemma 19.5]{EKM08}, at least one of the residue forms modulo $(\ovl{\pi})$ must be isotropic over $\kappa$.  
This shows (a), and by \Cref{unimodulargeneral}, this implies (b).
\end{proof}

Although we do not need it in the sequel of this article, we observe that the statement of \Cref{residueisotropic} can be made more symmetric when the regular local ring $R$ is additionally henselian.

\begin{prop}\label{two-dimhenselianresidue}
Let $R$ be a $2$-dimensional regular local ring with maximal ideal $\mf m$, fraction field $F$ and residue field $\kappa$ and $2\in\mg{R}$.
Assume that $R$ is henselian. 
Let $(\pi, \delta)$ be a parameter system of $R$ and let $a_1, \ldots, a_r \in  R$ be square-free and such that $\supp_R(a_1\cdots a_r) \subseteq \supp_R(\pi\delta)$.
Let $F^{(\delta)}$ and $F^{(\pi)}$ denote the completions of $R$ with respect to the height-$1$ prime ideals $(\delta)$ and $(\pi)$, respectively.
Then the following are equivalent:
\begin{enumerate}[$(i)$]
\item $\la a_1, \ldots, a_r\ra_F$ is isotropic.
\item $\la a_1, \ldots, a_r\ra_{F^{(\delta)}}$ is isotropic.
\item  $\la a_1, \ldots, a_r\ra_{F^{(\pi)}}$ is isotropic.
\end{enumerate}
\end{prop}

\begin{proof}
Since $F$ is contained in $F^{(\pi)}$ and $F^{(\delta)}$,  $(i)$  implies $(ii)$ and $(iii)$.

Assume that $(ii)$ holds.
By the hypothesis, we have 
\[a_1,\dots,a_r\in \left(\mg{R}\cup \pi\mg{R}\cup\delta\mg{R}\cup\pi\delta\mg{R}\right).\]
Hence, there exist unimodular diagonal forms $\varphi_0,\varphi_1,\varphi_2,\varphi_3 $ over $R$ such that, over $F$, we have  
$$\la a_1,\dots,a_r\ra_F\simeq  \varphi_0 \perp \pi\varphi_1 \perp \delta \varphi_2\perp \pi\delta \varphi_3. $$
It follows by \Cref{residueisotropic} that one of the residue forms of $\varphi_0,{\varphi}_1,{\varphi}_2,{\varphi}_3$ modulo $\mfm$ is isotropic  over $\kappa$. 
Since $R$ is henselian, we conclude that $\la a_1,\dots,a_r\ra_F$ is isotropic over $F$.
Hence we have shown that $(ii)$ implies $(i)$. 
An analogous argument  with switched roles of $\pi$ and $\delta$ yields that $(iii)$ implies $(i)$.
\end{proof}

\section{Function fields of surfaces}\label{S:FFS} 

In this section, we recall Lipman's results on resolution of singularities for surfaces and derive some statements tailored to our needs related to applications to quadratic forms over semiglobal fields in the next section.

Consider a domain $T$ with fraction field $K$ and a scheme $\mc X$ endowed with a  morphism $\eta: \mc X\to \Spec(T)$.
We call $\mc X \times_{\Spec(T)} \Spec (K)$ the \textit{generic fiber of $\eta$} and, if $T$ is local with residue field $\kappa$, then we call 
$\mc X \times_{\Spec(T)} \Spec (\kappa)$ the \textit{special fiber of $\eta$}.

\begin{lem}\label{proj-surj}
Let $T$ be a domain and $\mc{X}$ an integral scheme endowed with a projective morphism $\eta:\mc{X}\to \Spec(T)$ whose generic fiber is nonempty.
Then $T$ is contained in the function field of $\mc{X}$ and $\eta$ is surjective.
In particular, $\eta$ is constant if and only if $T$ is a field.
\end{lem}
\begin{proof}
Since the morphism $\eta$ is projective, its image is closed in $\Spec(T)$. 
By hypothesis, the image of $\eta$ contains the generic point of $\Spec(T)$. In particular, the generic point of $\mc{X}$ maps to the generic point of $\Spec(T)$.
Hence $\eta$ induces an embedding of $T$ into the function field of $\mc{X}$.
Since $\Spec(T)$ equals the closure of its generic point, we obtain that $\eta$ is surjective.
In particular, $\eta$ is constant if and only if $\Spec(T)$ consists of a single point, which holds if and only if $T$ is a field.
\end{proof}

\begin{lem}\label{dimension 2 argument}
Let $T$ be a Dedekind domain that is not a field.
 Let $\mc X$ be an integral separated noetherian scheme endowed with a projective morphism $\mc X\to \Spec(T)$ whose generic fiber is integral of dimension one.
Then $\dim(\mathcal{X})=2$. 
\end{lem}

\begin{proof}
We denote by $\eta:\mc{X}\to \Spec(T)$ the given morphism and by $\mc{X}_0$ its generic fiber.
By the hypothesis, $\dim \mc{X}_0=\dim(\Spec(T))=1$, and 
by \Cref{proj-surj}, $\eta$ is not constant.
As $\Spec(T)$ is a Dedekind scheme, we obtain by \cite[Corollary 4.3.10]{Liu12} that $\eta$ is flat.
We conclude by \cite[Theorem 4.3.12 \& Cor. 4.3.14]{Liu12} that
$\dim \mc{X} = \dim (\Spec(T)) + \dim \mc{X}_0=2$. 
\end{proof}

By a \textit{surface} we mean an integral separated noetherian scheme of dimension~$2$.
We will consider surfaces $\mc{X}$ endowed with a projective morphism $\mc{X}\to\Spec(T)$ where $T$ is a normal noetherian domain of dimension $1$ or $2$.

A surface $\mc{X}$ is \textit{regular} if, for every point $x \in \mc X$, the local ring $\mc{O}_{\mc X, x}$ is regular.
Let $\mc{X}$ be a regular surface. 
For $i\in\{0,1,2\}$, we denote by $\mc{X}^{(i)}$ the set of codimension-$i$ points in $\mc{X}$.
Denoting by $F$ the function field of $\mc{X}$, every $\xi\in\mc{X}^{(1)}$ determines a $\zz$-valuation $v_\xi$ on $F$ whose valuation ring is given by
$\mc{O}_{\mc{X},\xi}$.

Let $D$ be a (Cartier) divisor on $\mc{X}$.
For $\xi\in\mc{X}^{(1)}$, we set $v_\xi(D)=v_\xi(f)$ for any $f\in \mg{F}$ that locally defines $D$ at $\xi$.
The \textit{support of $D$ on $\mc{X}$} is the set 
$$\supp(D)= \{\xi \in \mc{X}^{(1)} \mid v_\xi(D)\neq 0\}.$$
We call the divisor $D$ \textit{effective} if $v_\xi(D)\geq 0$ for all $\xi\in\mc{X}^{(1)}$, and we call it \textit{reduced} if $v_\xi(D)\in\{-1,0,1\}$ for all $\xi\in\mc{X}^{(1)}$.
Any finite subset of $\mc{X}^{(1)}$ is equal to the support of a unique reduced effective divisor of $\mc{X}$.
Hence, effective divisors on $\mc{X}$ are in one-to-one correspondence to closed subschemes of $\mc{X}$ which are locally defined by nonzero principal ideals.
We also define the \textit{support} of a proper closed subset $V$ of $\mc{X}$ as 
$\supp(V)=\{\xi \in \mc{X}^{(1)} \mid \xi \in V\}.$
This is the (finite) set of generic points of the irreducible components of
  $V$ of codimension one in $\mc{X}$.
For $x\in\mc{X}^{(2)}$, we say that a divisor $D$ on $\mc{X}$ has
\textit{normal crossing at $x$} if  the $2$-dimensional local ring $\mathcal{O}_{\mc X,x}$ has a parameter system $(\pi,\delta)$ such that $D$ is locally defined by 
$\delta^l\pi^m$ for certain $l,m\in\zz$.
We say that \textit{$D$ has only normal crossings on $\mc{X}$} if $D$ has normal crossing at every point $x\in\mc{X}^{(2)}$.

Given an arbitrary surface $\mc{X}$, a \textit{regular model of $\mc{X}$} is a regular surface $\mc{X}'$ endowed with a proper birational morphism $\eta: \mc X' \to \mc X$, and we call it \textit{projective} if $\eta$ is projective. 
Given a regular scheme $\mc{X}$ and an effective divisor $D$ on $\mc{X}$, considered as a closed subscheme $D \hookrightarrow \mc{X}$, and $\eta: \mc{X}' \to \mc{X}$ a regular model of $\mc{X}$, we denote by
$\eta^{-1}(D)= D \times_{\mc{X}} \mc{X}'$ its scheme theoretic preimage, which corresponds to an effective divisor on $\mc{X}'$, by \cite[Tag 02OO]{Stacks}.

\smallskip

We now turn to the problem whether a given surface has a regular proper model.
J.~Lipman showed that, when the surface is excellent, such a model can be constructed as a composition of normalizations and blow-ups.
We refer to \cite[Subsection 4.1.2]{Liu12} for the notion of normalization, \cite[Section 8.1]{Liu12} for the construction of blow-up in a closed subscheme, and  \cite[Section 8.2]{Liu12} for the definition of excellent rings and schemes. We emphasize that excellent rings and schemes are noetherian by definition.
We give some relevant examples.

\begin{exs}\label{EX:excellent}
\begin{enumerate}[$(1)$]
\item Any complete noetherian local ring 
is excellent; see \cite[Theorem 8.2.39]{Liu12}.
\item Any discrete valuation ring of characteristic zero is excellent; see \cite[Cor.~8.2.40]{Liu12}.
\item Any projective scheme over an excellent scheme is excellent; see \cite[Theorem 8.2.39~(c)]{Liu12}.
\end{enumerate}
\end{exs}

\begin{prop}\label{blowup-normalization-projective}
Let $T$ be a noetherian ring and $\mc X$ an excellent scheme which is
projective over $\Spec(T)$. Let
$\eta:\mc X'\longrightarrow\mc X$
be a morphism which is a composition of normalizations and blow-ups. Then $\eta$ is projective and birational, and $\mc X'$ is excellent.
\end{prop}
\begin{proof}
Since blow-ups and normalizations are  birational morphisms, their composition $\eta$ is also birational.
By \cite[Cor.~3.3.32]{Liu12}, compositions of projective morphisms are projective.
Since $\mc X$ is assumed to be excellent, it therefore suffices to prove the statement for the isolated steps where
$\eta$ is either a normalization or a blow-up.

Suppose first that $\eta$ is the normalization of $\mc X$. 
Since $\mc X$ is excellent, \cite[Theorem~8.2.39~(d)]{Liu12} yields that $\eta$ is finite and 
\cite[Theorem~8.2.39~(c)]{Liu12}
that $\mc X'$ is excellent. 
By \cite[Tag~0B3I]{Stacks}, since $\eta$ is finite, it is projective.  

Suppose now that $\eta$ is the blow-up of $\mc X$ in a closed subscheme.
Then $\eta$ is projective, by \cite[Prop.~8.1.22]{Liu12}. 
In particular $\eta$ is locally of finite type, and hence,
by \cite[Theorem~8.2.39~(c)]{Liu12}, we obtain that $\mc{X}'$ is excellent.
\end{proof}

\begin{thm}[Lipman]\label{Lipman} 
Let $\mc X$ be an excellent surface. 
Then there exists a regular proper model  $\eta: \mc X' \longrightarrow \mc X$. 
If $\mc X$ is projective over $\Spec(T)$ for a noetherian ring $T$, then $\eta$
 may be chosen to be projective.
\end{thm}

\begin{proof}
Since $\mc X$ is an excellent surface, by \cite[Theorem]{Lip78}, there exists a regular proper model
$\eta:\mc{X}'\longrightarrow\mc{X}$, 
and by \cite[Tag 0BGP]{Stacks}, one may choose $\eta$ as a composition of blow-ups and normalizations. (This is also what the proof in \cite{Lip78}  explicitly shows.)
The final part then follows by \Cref{blowup-normalization-projective}.
\end{proof}

\begin{thm}[Lipman]\label{Lipman-embedded} 
Let $\mc X$ be a regular excellent surface and let $D$ be an effective divisor
on $\mc X$. Then there exists a regular proper model
$\eta:\mc X'\longrightarrow\mc X$
such that $\eta^{-1}(D)$ has only normal crossings on $\mc X'$.
If $\mc X$ is projective over $\Spec(T)$ for a noetherian ring $T$, then $\eta$
 may be chosen to be projective.
\end{thm}

\begin{proof} 
Since closed immersions are of finite type, for any irreducible component $Y$ of $D$, considered as a closed reduced subscheme of $\mc{X}$, we obtain by \cite[Theorem 8.2.39~(c)]{Liu12}  that $Y$ is excellent and thus by \cite[Theorem 8.2.39~(d)]{Liu12}  that the normalization of $Y$ is finite over $Y$.
Therefore, by \cite[Tag~0BIC]{Stacks}, there exists a composition of blow-ups at closed
points 
$\eta:\mc X'\longrightarrow\mc X$
such that $\eta^{-1}(D)$ has only normal
crossings on $\mc X'$. Moreover, since $\mc{X}$ is regular, we have by \cite[Theorem 8.1.29]{Liu12} that $\mc{X}'$ is regular.

If $\mc{X}$ is projective over $\Spec(T)$ for a noetherian ring $T$, then the blow-ups are projective, by \cite[Prop.~8.1.22]{Liu12}, and it follows inductively by using
\cite[Cor.~3.3.32]{Liu12} that $\eta$ is projective.
\end{proof}

We now show how to naturally realize the field extensions appearing in the definition of a semiglobal field as function fields of regular surfaces.

\begin{prop}\label{algfunctionfield} Let $d\in\{1,2\}$. Let $T$ be an excellent $d$-dimensional normal local domain and $K$ its fraction field. Let $F/K$ be a finitely generated field extension of transcendence degree $2-d$.
Then $F$ is the function field of a regular surface $\mc X$ admitting a surjective projective morphism $\mc X \rightarrow \Spec (T)$ such that all components of its special fiber have codimension $1$.
\end{prop}

\begin{proof} 
We first distinguish cases according to $d$ to construct a surface $\mc{Y}$ together with a projective and surjective morphism $\eta_Y:\mc{Y}\to \Spec(T)$.

Suppose first $d=1$. By \cite[Prop.~7.3.13]{Liu12} we have that $F$ is the function field of a regular projective curve $C$ over $K$. 
We choose $m\in\nat$ and homogeneous polynomials $f_1,\ldots,f_m \in T[X_0,\ldots,X_n]$ which together define $C$ inside $\mathbb{P}^n_K$. 
Let $\mc{Z}=\mathsf{Proj} (T[X_0,\ldots,X_n]/(f_1,\ldots,f_m))$.
Then $\mc{Z}$ is a projective noetherian scheme over $T$ with generic fiber $C$. 
The Zariski closure of $C$ in $\mc{Z}$ yields a closed subscheme $\mc{Y}$ of $\mc{Z}$ endowed with its induced reduced structure. 
It follows that $\mc{Y}$ is noetherian and integral. 
By \cite[Cor.~3.3.23]{Liu12}, the corresponding closed immersion $\mc{Y}\to\mc{Z}$ is projective.
Hence, its composition with the canonical structure morphism $\mc{Z}\to \Spec(T)$ yields a projective morphism $\eta_{\mc{Y}}:\mc{Y} \to  \Spec(T)$. 
In particular, $\eta_{\mc{Y}}$ is closed, and since its generic fiber is nonempty, \Cref{proj-surj} yields that it is surjective and nonconstant. 
Thus $\mc{Y}$ has dimension two, by \Cref{dimension 2 argument}. 
The affine scheme $\Spec(T)$ and the projective morphism $\eta_{\mc{Y}}$ are both separated, by  \cite[Prop. 3.3.4]{Liu12} and \cite[Cor.~3.3.10]{Liu12}, respectively. It follows by \cite[Prop. 3.3.9 (b)]{Liu12} that $\mc{Y}$ is separated. Hence $\mc{Y}$ is a surface.

Consider now the case where $d=2$.
Then $F/K$ is a finite field extension. 
We choose $m\in\nat$ and $x_1,\dots,x_m\in F$ integral over $T$ such that $F=K(x_1,\dots,x_m)$, and we set $T'=T[x_1,\dots,x_m]$ and $\mc{Y}=\Spec(T')$.
Since $T$ is noetherian, so is $T'$.
As an affine scheme, $\mc{Y}$ is separated  by \cite[Prop.~3.3.4]{Liu12}.
The natural morphism $\eta_Y:\mc{Y}\to\Spec(T)$ is finite, and in particular projective, according to \cite[Tag 0B3I]{Stacks}. 
Since $T'$ is integral over $T$, it follows by \cite[Prop.~4.15 \& Cor.~4.18]{Eisenbud} that $\eta_Y$ is surjective and $\dim(\mc{Y})=\dim(T')=\dim(T)=2$.
Hence $\mc{Y}$ is a surface.

In either case ($d=1$ or $d=2$), we have that $F/K$ is the function field of a surface $\mc{Y}$ endowed with a surjective projective morphism $\eta_{\mc{Y}}:\mc{Y}\to\Spec(T)$. 

Since $T$ is excellent and $\eta_{\mc{Y}}$ is in particular of finite type, we obtain by \cite[Theorem 8.2.39~(c)]{Liu12} that $\mc{Y}$ is excellent. 
Hence, by \Cref{Lipman}, there exists a regular projective model $\eta: \mc X \to \mc{Y}$. 
We may further choose this model in such way that all components of the special fiber of the composition $\eta_{\mc{Y}}\circ \eta:\mc{X}\to\Spec(T)$ have codimension $1$.
Indeed, if for a first choice of the model this fails, then we may compose $\eta$ with the blow-up in $\mc{X}$ in the closed subscheme given by the isolated closed points in the special fiber of $\eta_{\mc{Y}}\circ \eta$ to achieve a regular projective model with the desired property.

We now consider the morphism $\eta_{\mc{X}}=\eta_{\mc{Y}}\circ\eta:\mc{X}\to \Spec(T)$.
Since $\eta$ is birational, its image contains the generic point of $\mc{Y}$.
Since $\eta_{\mc{Y}}$ is surjective, it follows that the generic fiber of $\eta_{\mc{X}}$ is nonempty.
Being a composition of projective morphisms, $\eta_{\mc{X}}$ is projective, by \cite[Cor.~3.3.32~(b)]{Liu12}, and we conclude by \Cref{proj-surj} that it is also surjective. 
\end{proof}

Having constructed suitable regular projective models, we conclude this section by analyzing the local rings and residue fields of specific points on these surfaces. These geometric properties will serve in \Cref{S:RMQF} as the key bridge to the local algebraic results of \Cref{S:prelim}.

\begin{lem}\label{finite field extension}
Let $T$ be a noetherian $2$-dimensional normal local domain and let $E$ be its fraction field.
Let $\mc{X}$ be a surface endowed with  a  surjective projective morphism $\eta:\mc X \rightarrow \Spec(T)$. Let $F$ be the function field of $\mc{X}$.
Then $F/E$ is a finite field extension.
\end{lem}

\begin{proof}
It suffices to show that the generic fiber of $\eta$ has dimension zero. Let $\mc{X}_s$ denote the special fiber of $\eta$.
Consider the open subscheme $\mc Y= \mc{X} \setminus \mc{X}_s$ of $\mc{X}$. 
Let $\mfm$ be the maximal ideal of $T$.
The open subscheme $S=\Spec(T) \setminus \{\mfm\}$ of $\Spec(T)$ is a Dedekind scheme, and $\eta$ restricts to a surjective  morphism $\eta_{\mc{Y}}:\mathcal{Y} \to S$.
By \cite[Cor.~4.3.10]{Liu12}, it follows that $\eta_{\mc{Y}}:\mc{Y}\to S$ is flat. 
Since $\eta$ is projective, and in particular closed, $\mc Y$ does not contain any closed points of $\mc{X}$.
Since $S$ is open in $\Spec(T)$ and $\mc{Y}=\eta^{-1}(S)$, we obtain that $\eta_{\mc{Y}}:\mc{Y}\to S$ is projective.

Consider an arbitrary closed point $y$ of $\mc{Y}$ and set $\mfp=\eta_{\mc{Y}}(y)$.
Then $\mfp$ is closed in $S$.
Hence $\mfp$ is a height-one prime ideal of $T$ and we have $\dim \mc{O}_{\mc Y,y} =1= \dim T_\mfp$.
Since $\eta_{\mc{Y}}$ is flat, we obtain by \cite[Theorem 4.3.12]{Liu12}  that $\dim \mathcal{O}_{\mc Y_\mfp,y}=\dim \mc{O}_{\mc Y,y} - \dim T_{\mfp} =0$, where $\mc Y_\mfp$ is the fiber of $\eta_{\mc Y}$ over $\mfp$.
From this it follows by \cite[Prop.~4.4.16]{Liu12} that the generic fiber of $\eta_{\mc{Y}}:\mc{Y}\to S$ has dimension zero. Since $F$ is the function field of $\mc Y$ and $E$ the function field of $S$, we conclude that the field extension $F/E$ is finite.
\end{proof}

\begin{lem} \label{unipointcurves}
Let $T$ be a henselian normal local noetherian domain of dimension $1$ or $2$. 
Let $\mc{X}$ be a surface endowed with a surjective projective morphism $\eta:\mc X \rightarrow \Spec(T)$.
Consider a point $x\in \mc X^{(1)}$ which is not contained in the special fiber of $\eta$.
Then the field $\kappa_x$ carries a henselian $\zz$-valuation. 
\end{lem}

\begin{proof}
We show that $\kappa_x$ is the fraction field of a henselian discrete valuation ring.
Let $E$ denote the fraction field of $T$. 

We first consider the case where $\dim T=1$. 
Then $x$ lies in the generic fiber of $\eta$. Since all points of $\mc{X}^{(2)}$ are contained in the special fiber of $\eta$, it follows that $\{x\}$ is closed in the generic fiber.
We obtain that $E\subseteq \mc{O}_{\mc{X},x}$, whereby $\kappa_x$ is a finite  extension of $E$.
As $\dim(T)=1$, $T$ is a henselian discrete valuation ring, and hence so is its integral closure in $\kappa_x$.

Assume now that $\dim T=2$.
Let $F$ denote the function field of $\mc{X}$.
By \Cref{finite field extension},  $F/E$ is a finite field extension. 
Then
$\mathfrak{p}= \eta(x)$ is a height-one prime ideal of $T$,  $\kappa_x/ \kappa_\mathfrak{p}$ is a finite field extension and $T/\mf p$ is a henselian local domain of Krull dimension one with fraction field $\kappa_{\mf p}$. 
Let $B$ be the integral closure of $T/\mf p$ in $\kappa_x$. 
By \Cref{KA}, $B$ is a henselian discrete valuation ring of $\kappa_x$.
\end{proof}

\begin{lem}\label{prelem}
Let $T$ be a noetherian domain of dimension $1$ or $2$ with $2\in \mg T$.
Let $\mc X$ be a regular surface endowed with a projective surjective morphism $\mc X \rightarrow \Spec(T)$. Let $F$ denote the function field of $\mc{X}$.
Let $v$ be a  valuation on $F$ such that $T\subseteq \mc{O}_v$.
Then there exists $x\in \mc X$ such that $\mc O_{\mc X,x} \subseteq \mc O_v$ with $\mf m_v \cap \mc O_{\mc X,x} = \mf m_{\mc X,x}$.
\end{lem}

\begin{proof}
Let $\xi$ denote the generic point of $\mc X$, so that $\mc O_{\mc X,\xi} = F$.
Since the given morphism $\eta\colon\mc{X}\to\Spec(T)$ is surjective,
$\eta(\xi)$ is the generic point of $\Spec(T)$, so $T$ embeds into $F$ via $\eta$.
As $\eta$ is projective, it is proper, by \cite[Tag~01WC]{Stacks}.
Hence, by the valuative criterion of properness \cite[Tag~0BX5]{Stacks},
there exists a morphism $\psi\colon\Spec(\mc O_v) \to \mc X$ making the diagram
\[
\begin{tikzcd}
\Spec(F) \arrow[r] \arrow[d] & \mc X \arrow[d, "\eta"] \\
\Spec(\mc O_v) \arrow[r] \arrow[ur, dashed, "\psi"] & \Spec(T)
\end{tikzcd}
\]
commute, where the upper horizontal arrow is the canonical morphism with
image $\{\xi\}$.
Let $x \in \mc X$ be the image under $\psi$ of the closed point of
$\Spec(\mc O_v)$. Then $\psi$ induces on stalks a local homomorphism
$\varphi\colon \mc O_{\mc X,x} \to \mc O_v$.
Since the upper triangle commutes, $\psi$ maps the generic point of
$\Spec(\mc O_v)$ to $\xi$, and composing $\varphi$ with the embedding
$\mc O_v\to F$ gives the canonical map
$\mc O_{\mc X,x} \to \mc O_{\mc X,\xi} = F$, which is injective
as $\mc X$ is integral.
Hence $\varphi$ is injective, whereby $\mc O_{\mc X,x} \subseteq \mc O_v$.
Since $\varphi$ is local, we obtain that
$\mf m_v \cap \mc O_{\mc X,x} = \mf m_{\mc X,x}$.
\end{proof}

\section{Regular models adapted to a quadratic form}\label{S:RMQF}

Given a quadratic form over a semiglobal field,
we use the flexibility in the choice of projective models established in Section 3 to fit the local hypotheses in the local isotropy criteria of Section 2.

For a surface $\mc{X}$ endowed with a morphism $\eta:\mc{X}\to \Spec(T)$ for a local domain $T$, we denote by $\mc{X}_s$ the special fiber of $\eta$.

\begin{prop}\label{refinement}
Let $T$ be a normal noetherian henselian local domain of dimension $1$ or $2$ with $2\in \mg T$.
Let $\mc{X}$ be a regular surface endowed with a projective and surjective morphism $\mc{X}\to \Spec(T)$ and such that all the irreducible components of $\mc{X}_s$ are $1$-dimensional.
Let $F$ denote the function field of $\mc{X}$.
Let $n\in\nat$ and $a_1,\dots,a_n\in\mg{F}$ be such that the reduced effective divisor $D$ with support $\supp_{\mc X}(a_1) \cup \ldots \cup \supp_{\mc X}(a_n) \cup \supp_{\mc X}(\mc X_s)$ has normal crossings on $\mc X$. 
If $\la a_1, \ldots, a_n \ra_{F^{v_y}}$ is isotropic for every $y \in \mc{X}^{(1)} \cap \mc{X}_s$,
then $\la a_1, \ldots, a_n \ra_{F^v}$ is isotropic for every rank-$1$ valuation $v$ on $F$ with $T\subseteq \mc{O}_v$.
\end{prop}

\begin{proof}
Let $\vf=\la a_1, \ldots, a_n \ra$.
Suppose that there exists a rank-$1$ valuation $v$ on $F$ with $T\subseteq \mc{O}_v$ such that $\vf_{F^v}$ is anisotropic. 
By \Cref{prelem}, there exists $x\in \mc X$ such that $\mc O_{\mc X,x} \subseteq \mc O_v$ with $\mf m_v \cap \mc O_{\mc X,x} = \mf m_{\mc X,x}$.
If $x\in\mc{X}^{(1)}\cap \mc{X}_s$, then we may take $y=x$.

Suppose next that $x \in \mc X^{(2)}$. Let $\delta$ be a local equation of an irreducible component of the special fiber at the closed point $x$. Then by \Cref{residueisotropic}, we get that $\varphi$ is anisotropic over $F^{(\delta)}$.

Suppose now that $x \in \mc X^{(1)}\setminus\mc{X}_s$. Then $\mc{O}_{\mc{X},x}$ is a discrete valuation ring, and  hence it follows that $\mc{O}_{\mc{X},x}=\mc{O}_v$. 
Let $P$ be a closed point on the prime divisor $\ovl{\{x\}}$ on $\mc{X}$. Then $P$ is a closed point of $\mc{X}$ and hence $P\in\mc{X}_s$, since projective morphisms are closed maps.  
Let $\delta$ be a local equation of an irreducible component of the special fiber at $P$. 
By \Cref{unipointcurves}, $\kappa_x$ carries a henselian discrete valuation.
Since $\kappa_x=Fv$, it follows by \Cref{residueisotropic} that $ \varphi$ is anisotropic over $F^{(\delta)}$.

In either of the last two cases, for a sufficiently small affine neighbourhood of the closed point $x$ or $P$,  $\delta$ defines a point $y\in\mc{X}^{(1)}\cap\mc{X}_s$ such that $\vf_{F^{v_y}}$ is anisotropic.
\end{proof}

We are ready to prove our main result.

\begin{thm}\label{1-2-dim}
Let $F$ be a field and $T$ a normal excellent henselian local subring of $F$ with $2\in\mg{T}$ and such that one of the following holds:
\begin{enumerate}[$(i)$]
    \item $T$ is $1$-dimensional and $F$ is a function field in one variable over the fraction field of $T$.
    \item $T$ is $2$-dimensional and $F$ is its fraction field.
\end{enumerate}
Let $\vf$ be a quadratic form over $F$. Assume that $\vf$ is isotropic over $F^v$ for every $v \in \Omega^\circ_T(F)$.
Then $\vf$ is isotropic over $F^v$ for every rank-$1$ valuation $v$ on $F$  such that $T\subseteq \mathcal{O}_v$.
\end{thm}
\begin{proof}
    
By \Cref{algfunctionfield}, we have that $F$ is the function field of a regular surface $\mc X$ admitting a surjective projective morphism $\eta:\mc X \rightarrow \Spec (T)$ such that all components of $\mc{X}_s$ are of codimension one in $\mc{X}$.
We may assume that $\vf= \la a_1,\ldots,a_n\ra$ for some $n\in \mathbb{N}$ and $a_1, \ldots, a_n \in \mg{F}$.  
Consider the reduced effective divisor $D$ with
\[\supp_{\mc{X}}(D)=\supp_{\mc X}(a_1) \cup \ldots \cup \supp_{\mc X}(a_n) \cup \supp_{\mc X}(\mc X_s).\] 
By \Cref{Lipman-embedded}, there exists a regular projective model $\phi: \mc{X}' \to \mc{X}$ such that $\phi^{-1}(D)$ has only normal crossings on $\mc{X}'$. 
Let $\eta':\mc{X}'\to\Spec(T)$ be the composition of $\eta$ with $\phi$.
Being a composition of projective morphisms, $\eta'$ is projective, by \cite[Cor.~3.3.32~(b)]{Liu12}.
Since $\phi$ is birational, the generic fiber of $\eta'$ is nonempty.
Hence $\eta'$ is surjective, by \Cref{proj-surj}.
Let $D'$ be the reduced effective divisor on $\mc{X}'$ with \[\supp_{\mc{X}'}(D')=\supp_{\mc X'}(a_1) \cup \ldots \cup \supp_{\mc X'}(a_n) \cup \supp_{\mc X'}(\mc X_s').\]
Since $\supp(D')\subseteq\supp(\eta^{-1}(D))$, the divisor $D'$ has only normal crossings.
The statement now follows by \Cref{refinement} applied to $\mc{X}'$ and $D'$.
\end{proof}

We apply our results to retrieve the local-local principle for quadratic forms of dimension at least $3$.

\begin{cor}\label{2-foldlgpi}
Let $F$ be a field and $T$ a normal excellent henselian local subring of $F$ with $2\in\mg{T}$ and such that one of the following holds:
\begin{enumerate}[$(i)$]
\item $T$ is $1$-dimensional and $F$ is a function field in one variable over the field of fractions of $T$.
\item $T$ is $2$-dimensional, complete, either regular or equicharacteristic, and $F$ is the fraction field of $T$.
\end{enumerate}
Let $\vf$ be a quadratic form of dimension at least $3$ over $F$. Then $\vf$ is isotropic over $F$ if and only if it is isotropic over $F^v$ for every $v\in\Omega_T^\circ(F)$.
\end{cor}

\begin{proof}
In the case $(ii)$, the statement follows from \cite[Cor.~4.7]{HHK15} by applying \Cref{1-2-dim}.

Assume now that we are in case $(i)$.
If $\vf$ is isotropic over $F$, then also over $F^v$ for any valuation $v$ on $F$. 
Assume now that $\vf$ is anisotropic over $F$.
We now apply Mehmeti's local-global principle from \cite{Meh19} in the  extended version (relaxing the hypothesis on the base field from complete to henselian) reached in \cite[Theorem 4.4]{BDGMZ26}.
We obtain that $\varphi$ is anisotropic over $F^v$ for some $v\in \Omega(F)$ with $T \subseteq \mc{O}_v$. 
Then \Cref{1-2-dim} shows the existence of a valuation $w\in\Omega_T^\circ(F)$ such that $\vf$ is anisotropic over $F^w$.
\end{proof}

As a special case of \Cref{1-2-dim}, we further obtain the following local-local principle for sums of squares.

\begin{cor}\label{C:local-squares}
Under the hypothesis of \Cref{1-2-dim}, we have $$\mg{F}\cap\bigcap_{v\in \Omega(F)}\sq{(F^v)}=\mg{F}\cap\bigcap_{v\in\Omega_T^\circ(F)}\sq{(F^v)}\,.$$
\end{cor}

\begin{proof}
    For $a\in\mg{F}$ and a field extension $L/F$, the binary quadratic form $X_1^2-aX_2^2$ is isotropic over $L$ if and only if $a\in\sq{L}$. Hence the statement follows from \Cref{1-2-dim}.
\end{proof}

In a separate work \cite{BGMF26}, the group described in \Cref{C:local-squares} is studied further and identified as the Kaplansky radical of $F$.

\subsection*{Acknowledgements} 
We acknowledge support from the \textit{Bijzonder Onderzoeksfonds, Universiteit Antwerpen} (BOF-DOCPRO4, 2533) and the \textit{Universidad de Santiago de Chile} (Proyecto DICYT,
c\'odigo 042432G).

\section*{Declarations}

\subsection*{Data availability}
The  article describes entirely theoretical research.
All new data supporting the findings presented in this article are included in the manuscript.
Data sharing is not applicable to this article.

\subsection*{Conflict of interest} The authors declare that there is no conflict of interest.

\subsection*{AI usage statement}
The authors declare that they used AI tools (\emph{ChatGPT}, \emph{Google Gemini}, \emph{Claude}) during the finalization process for this article to identify inconsistencies and lack of clarity.

\bibliographystyle{amsalpha}

\addcontentsline{toc}{chapter}{Bibliography}

\begin{thebibliography}{99}

\bibitem{BDGMZ26} K.J.~Becher, N.~Daans, D.~Grimm, G.~ Manzano-Flores, M.~Zaninelli. The Pythagoras number of a rational function field in two variables. 
Preprint (2026),\\
\url{https://arxiv.org/abs/2302.11425}

\bibitem{BGVG14}
K.J.~Becher, D.~Grimm, J.~Van Geel. Sums of squares in algebraic function fields over a complete discretely valued field. \textit{Pacific J.~Math.} {267} (2014), 257--276.

\bibitem{BGMF26}
K.J.~Becher, D.~Grimm, G.~Manzano-Flores.
The Kaplansky radical of a function field in one variable. 
Preprint (2026).

\bibitem{BOU-AC5to7}
N.~Bourbaki.
\newblock \textit{\'El\'ements de math\'ematique. Alg\`ebre Commutative. Chapitres 5 \`a
  7. Reprint of the 1975 original.}
\newblock {Berlin: Springer}, 2006.

\bibitem{CTPS12}
J.-L.~Colliot-Thélène, R.~Parimala,  V.~Suresh.
Patching and local-global principles for homogeneous spaces over function fields of $p$-adic curves.
\textit{Comment. Math. Helv.} 87 (2012), 1011--1033.

\bibitem{Eisenbud}
D. Eisenbud. \textit{Commutative algebra: with a view toward algebraic geometry.}
Springer Graduate Texts in Mathematics, Volume 150, Springer, 2013.

\bibitem{EKM08}
R. Elman, N. Karpenko, A. Merkurjev. 
\textit{The algebraic and geometric theory of quadratic forms}. 
Amer.~Math.~Soc.~Coll.~Pub., vol. {56}, Amer. Math.~Soc., Providence, RI, 2008.

\bibitem{ELP73}
R.~Elman, T.Y.~Lam, A.~Prestel.
On some Hasse Principles over formally real fields.
\textit{Math. Z.} 134  (1973), 291--301.

\bibitem{EP05} 
A.J.~Engler, A.~Prestel. \textit{Valued Fields}. Springer, 2005.

\bibitem{Gup21}
P.~Gupta.
Four-dimensional quadratic forms over $C(\!(t)\!)(X)$.
\textit{Arch. Math. (Basel)} 117 (2021), 369--374.

\bibitem{HH10}
D. Harbater, J. Hartmann.
Patching over fields.
\textit{Israel J. Math.}, {176} (2010), 61--107.

\bibitem{HHK09}
D. Harbater, J. Hartmann, D. Krashen.
Applications of patching to quadratic forms and central simple algebras. 
\textit{Invent. Math.} 178 (2009), 231--263.

\bibitem{HHK15}
D.~Harbater, J.~Hartmann,  D.~Krashen. Refinements to patching and applications to field invariants. \textit{Int. Math. Res. Not.}
  2015 (2015), 10399--10450.

\bibitem{Kne76}
M.~Knebusch. On algebraic curves over real closed fields, II.
\textit{Math.~Z.} 151  (1976), 189--205.

\bibitem{Lip78}
J.~Lipman. Desingularization of two-dimensional schemes. \textit{Ann. Math.} 107 (1978), 151--207.

\bibitem{Liu12}
Q.~Liu. \textit{Algebraic geometry and arithmetic curves}.
Oxford Graduate Texts in 
Mathematics, Oxford University Press, 2012.

\bibitem{Meh19}
V. Mehmeti.
\newblock{Patching over Berkovich curves and quadratic forms.}
\textit{Comp. Math.} 155 (2019), 2399--2438.

\bibitem{Nag75}
M.~Nagata. \textit{Local rings}.
{Corr.~reprint of the orig.~ed. 1962}.
Robert E.~Krieger Publ. Co.~XIII, Huntington, New York, 1975.

\bibitem{Wit37}
E.~Witt. Theorie der quadratischen Formen in beliebigen Körpern. \textit{J. reine angew. Math.} 176  (1937),
31--44.

\bibitem{Stacks}
The Stacks Project Authors.
\textit{The Stacks Project}.
\url{https://stacks.math.columbia.edu}.

\end{thebibliography}

\end{document}